\documentclass[twocolumn]{autart}

\usepackage[utf8]{inputenc}
\usepackage{tabularx}
\usepackage{float}
\usepackage[caption=false]{subfig}
\usepackage{epsfig}
\usepackage{epstopdf}
\usepackage{bbm}
\usepackage{epsfig}
\usepackage{amsmath}
\usepackage{amssymb}
\usepackage{bm}
\usepackage{epstopdf}
\usepackage{color}

\usepackage{cite}
\usepackage{algorithm}
\usepackage{algorithmicx}
\usepackage{algpseudocode}

\usepackage{url}
\usepackage{enumitem}

\begin{document}

\begin{frontmatter}

\title{Data-Driven Predictive Control With Adaptive Disturbance Attenuation For Constrained Systems}
\author[NL]{Nan Li}\ead{nanli@auburn.edu},
\author[IK]{Ilya Kolmanovsky}\ead{ilya@umich.edu},
\author[HC]{Hong Chen}\ead{chenhong2019@tongji.edu.cn}

\address[NL]{Department of Aerospace Engineering, Auburn University, Auburn, AL 36849, USA}
\address[IK]{Department of Aerospace Engineering, University of Michigan, Ann Arbor, MI 48105, USA}
\address[HC]{College of Electronic and Information Engineering, Tongji University, Shanghai, China}
\begin{abstract}                          
In this paper, we propose a novel data-driven predictive control approach for systems subject to time-domain constraints. The~approach combines the strengths of $\mathcal{H}_{\infty}$ control for rejecting disturbances and MPC for handling constraints. In particular, the approach can dynamically adapt $\mathcal{H}_{\infty}$ disturbance attenuation performance depending on measured system state and forecasted disturbance level to satisfy constraints. We establish theoretical properties of the approach including robust guarantees of closed-loop stability, disturbance attenuation, constraint satisfaction under noisy data, as well as sufficient conditions for recursive feasibility, and illustrate the approach with a numerical example.
\end{abstract}
\end{frontmatter}

\section{Introduction}\label{sec:1}

In addition to stability, disturbance rejection and constraint satisfaction are important topics in control systems engineering, especially as modern systems often operate in complex and dynamic environments while increasingly stringent safety and performance requirements are imposed on them. The $\mathcal{H}_{\infty}$ control is an effective approach to addressing the former -- it aims to minimize the effect of disturbances on system outputs \cite{zhou1996robust}. Since introduced in the early 1980s, $\mathcal{H}_{\infty}$ control has been successfully implemented in aerospace, automotive, and many other sectors. However, time-domain constraints on state and control variables are not handled in conventional $\mathcal{H}_{\infty}$ control. Meanwhile, Model Predictive Control~(MPC) stands out from many other control approaches for its ability to explicitly and non-conservatively handle time-domain constraints \cite{maciejowski2002predictive}. Therefore, combing $\mathcal{H}_{\infty}$ control and MPC to achieve desired disturbance rejection properties while satisfying constraints is appealing and has been studied in, e.g., \cite{chen2004disturbance,chen2006moving,lee2008robust,orukpe2011towards,huang2014mixed,benallouch2014h,song2016mixed,song2017robust,zhang2018robust,shokrollahi2021robust}. A major strength of such a combination is the ability to dynamically adapt $\mathcal{H}_{\infty}$ disturbance attenuation performance (depending on measured/estimated system state and forecasted disturbance level) to satisfy constraints through solving an MPC optimization problem repeatedly in an online manner \cite{chen2004disturbance,chen2006moving}.

Conventional $\mathcal{H}_{\infty}$ control and MPC methods rely on parametric models of the system to be controlled. As engineered systems are becoming increasingly complex and cyber-physical, first-principle models are more difficult to obtain. Meanwhile, with the rapid advances in sensing, computation, and communication technologies, data is more readily available. This has spurred the development of data-driven methods for system modeling, analysis, and control. In particular, practitioners may favor an end-to-end solution that bypasses the intermediate steps of modeling and analysis and produces a controller with desired properties directly from measured data of system behavior. Therefore, it is beneficial to extend model-based methods for $\mathcal{H}_{\infty}$~control, MPC, and the aforementioned combination to their data-driven counterparts.

Here we provide a brief review of existing data-driven control methods in the literature: Reinforcement Learning~(RL) can be used to train optimal controllers from data~\cite{lewis2012reinforcement,recht2019tour}. However, the majority of RL methods optimize the average behavior of the closed-loop system when it is subject to disturbances with certain statistics and do not provide a worst-case robustness guarantee. Furthermore, they typically handle constraints through penalties, hence making constraints soft. Integrating data-driven uncertainty estimation/system identification with(in) the MPC framework can lead to reliable usage of data for improving MPC performance while maintaining certain robustness guarantees especially for satisfying constraints \cite{rosolia2017learning,hewing2020learning}. Along these lines, the Data-enabled Predictive Control (DeePC) is an emerging technique and has attracted increasing attention from researchers recently. The DeePC uses measured input-output data to create a non-parametric system model based on behavioral systems theory and uses the model to predict future trajectories~\cite{coulson2019data}. It has demonstrated superior performance in various applications~\cite{elokda2021data,huang2021decentralized,chinde2022data,li2023minimum}. However, DeePC entails higher computational cost than conventional MPC due to its high-dimensional, data-based non-parametric system model \cite{baros2022online,dai2022cloud}. Also, how to handle noisy data and provide certain robustness guarantees under noisy data remains to be an open question in DeePC \cite{berberich2020data,coulson2021distributionally,huang2023robust}. Furthermore, to the best of our knowledge, there has been no previous work integrating $\mathcal{H}_{\infty}$ control and MPC (or, designing an MPC that has an $\mathcal{H}_{\infty}$-type disturbance attenuation property) in the general data-driven MPC literature. An approach to synthesizing an $\mathcal{H}_{\infty}$ controller using noisy data based on a matrix S-lemma was introduced in \cite{van2020noisy}. The approach reduces a data-driven $\mathcal{H}_{\infty}$ control synthesis problem to a low-dimensional Linear Matrix Inequality (LMI) optimization problem, which is computationally tractable with state-of-the-art interior-point LMI solvers. However, time-domain constraints are not handled by the approach of \cite{van2020noisy}, and moving-horizon, MPC-type implementation of the approach was not considered in~\cite{van2020noisy}. More classical data-driven control methods also include self-tuning regulators \cite{aastrom1977theory} and iterative learning control \cite{bristow2006survey}. A comprehensive survey of data-driven control methods can be found in \cite{hou2013model}.

In this paper, we fill the gap in the literature by proposing a novel data-driven control approach that combines the strengths of $\mathcal{H}_{\infty}$ control for rejecting disturbances and MPC for handling constraints. Our approach can be viewed as a data-driven counterpart of the model-based moving-horizon $\mathcal{H}_{\infty}$ control approach of \cite{chen2004disturbance,chen2006moving} and enjoys similar properties including dynamic adaptation of $\mathcal{H}_{\infty}$ performance depending on measured system state and forecasted disturbance level for satisfying constraints. Specifically, the contributions include:
\begin{enumerate}
    \item Our approach is the first data-driven MPC method in the literature that focuses on $\mathcal{H}_{\infty}$-type disturbance attenuation for systems with time-domain constraints.
    \item We conduct a comprehensive analysis of the theoretical properties of our approach. The results include robust guarantees of closed-loop stability, disturbance attenuation, and constraint satisfaction under noisy data, as well as conditions for online problem recursive feasibility.
\end{enumerate}

The paper is organized as follows: We describe the problem treated in this paper including key assumptions and preliminaries in Section~\ref{sec:2}. We develop an approach to synthesizing an $\mathcal{H}_{\infty}$ controller that also enforces time-domain constraints for an unknown system using noisy trajectory data in Section~\ref{sec:3}. The development in this section has merit in its own right because it extends the data-driven $\mathcal{H}_{\infty}$ control synthesis approach of \cite{van2020noisy} (which does not handle constraints) to the constrained case. It is also an essential building block of the data-driven MPC algorithm developed in Section~\ref{sec:4}. Then, Section~\ref{sec:4} presents our proposed data-driven MPC algorithm and analyzes its theoretical properties. We illustrate the algorithm with a numerical example in Section~\ref{sec:5}. Finally, we conclude the paper in Section~\ref{sec:6}.

The notations used in this paper are mostly standard. We use $\mathbb{R}^n$ to denote the space of $n$-dimensional real vectors, $\mathbb{R}^{n \times m}$ the space of $n$-by-$m$ real matrices, and $\mathbb{N}$~the~set of natural numbers including zero. Given a vector $x \in \mathbb{R}^n$, we use $\|x\|$ to denote its Euclidean norm, i.e., $\|x\| = \sqrt{x^{\top} x}$. Given a matrix $M \in \mathbb{R}^{n \times m}$, its kernel is the subspace of all $x \in \mathbb{R}^m$ such that $Mx = 0$, i.e., $\text{ker}(M) = \{x \in \mathbb{R}^m: Mx = 0\}$. Given two symmetric matrices $M, N \in \mathbb{R}^{n \times n}$, $M \succ N$ means that $M - N$ is~positive definite, i.e., $x^{\top} (M - N) x > 0$ for all non-zero $x \in \mathbb{R}^n$, and $M \succeq N$ means that $M - N$ is positive semidefinite, i.e., $x^{\top} (M - N) x \ge 0$ for all $x \in \mathbb{R}^n$. Similarly, $M \prec N$ and $M \preceq N$ means that $M - N$ is negative definite and negative semidefinite, respectively. For an optimization problem, by ``solution'' we mean a feasible solution, i.e., a set of values for the decision variables that satisfies all constraints, and by ``optimal solution'' we mean a feasible solution where the objective function (almost) reaches its maximum (or minimum) value. Because the optimization problems appearing in this paper are all convex problems, an ``optimal solution'' is globally optimal.

\section{Problem Statement and Preliminaries}\label{sec:2}

We consider the control of dynamic systems that can be represented by the following linear time-invariant model:
\begin{subequations}\label{equ:1}
\begin{align}
x(t+1) &= A_{\text{o}} x(t) + B_{\text{o}} u(t) + w(t) \label{equ:1_1} \\
y_1(t) &= C_1 x(t) + D_1 u(t) \\
y_2(t) &= C_2 x(t) + D_2 u(t)
\end{align}
\end{subequations}
where $x(t) \in \mathbb{R}^n$ denotes the system state at the discrete time $t \in \mathbb{N}$, $u(t) \in \mathbb{R}^m$ denotes the control input, $w(t) \in \mathbb{R}^n$ denotes an unmeasured disturbance input, $y_1(t) \in \mathbb{R}^{p_1}$ and $y_2(t) \in \mathbb{R}^{p_2}$ are two outputs the roles of which are introduced below. The goal is to design a control algorithm that achieves the following three objectives:
\begin{enumerate}
    \item Guaranteeing closed-loop stability;
    \item Optimizing disturbance attenuation in terms of the $\mathcal{H}_{\infty}$ performance from $w$ to output $y_1$;
    \item Enforcing the following constraints on output $y_2$ at all times $t \in \mathbb{N}$:
    \begin{equation}\label{equ:2}
    y_{2v}(t) \le y_{2v,\max}, \quad v = 1,2,\dots,p_2
    \end{equation}
    where $y_{2v}(t)$ denotes the $v$th entry of $y_2(t)$ and $y_{2v,\max} \ge 0$ for all $v = 1,2,\dots,p_2$. Constraints in this form can represent state/output variable bounds, control input limits, etc.
\end{enumerate}

In addition to the standard assumptions of $(A_{\text{o}},B_{\text{o}})$ being stabilizable and $(C_1,A_{\text{o}})$ being detectable, we also assume that the system model $(A_{\text{o}},B_{\text{o}})$ is unknown and only trajectory data are available. This calls for a data-driven control approach.

Assume we have data $\{(x^+_j,x_j,u_j)\}_{j = 1}^{J}$, where $(x_j,u_j)$ denotes a pair of previous state and control input values, $x^+_j$ denotes the corresponding next state value, and the subscript $j$ indicates the $j$th data point. The data can be collected from a single or multiple trajectories. According to \eqref{equ:1_1}, we have
\begin{equation}\label{equ:3}
x^+_j = A_{\text{o}} x_j + B_{\text{o}} u_j + w_j
\end{equation}
where $w_j$ is the disturbance input value at the time of collecting the data point $(x^+_j,x_j,u_j)$. Organizing data with the following matrices:
\begin{subequations}\label{equ:4}
\begin{align}
X^+ &= \left[x^+_1,\, x^+_2,\, \cdots,\, x^+_J \right] \\
X &= \left[x_1,\, x_2,\, \cdots,\, x_J \right] \\
U &= \left[u_1,\, u_2,\, \cdots,\, u_J \right] \\
W &= \left[w_1,\, w_2,\, \cdots,\, w_J \right] \label{equ:4_4}
\end{align}
\end{subequations}
the relation \eqref{equ:3} implies 
\begin{equation}\label{equ:5}
W = X^+ - A_{\text{o}} X - B_{\text{o}} U
\end{equation}
Because the disturbance input $w_j$ is not measured, $W$~in~\eqref{equ:4_4} and \eqref{equ:5} is unknown.  We make the following assumption about the data:

{\bf Assumption~1.} The disturbance input values $w_1,w_2,$ $\dots,w_J$, collected in $W$, satisfy the following quadratic matrix inequality:
\begin{equation}\label{equ:6}
\begin{bmatrix}
I \\ W^{\top} 
\end{bmatrix}^{\top}\begin{bmatrix}
\Phi_{11} & \Phi_{12} \\ \Phi_{12}^{\top} & \Phi_{22}  
\end{bmatrix} \begin{bmatrix}
I \\ W^{\top} 
\end{bmatrix} \succeq 0
\end{equation}
where $\Phi_{11} = \Phi_{11}^{\top} \in \mathbb{R}^{n \times n}$, $\Phi_{12} \in \mathbb{R}^{n \times J}$, and $\Phi_{22} = \Phi_{22}^{\top} \prec 0$ are known matrices.

When $\Phi_{12} = 0$ and $\Phi_{22} = -I$, \eqref{equ:6} reduces to $\sum_{j = 1}^J w_j w_j^{\top} \preceq \Phi_{11}$, which has the interpretation that the total energy of the disturbance inputs in data is bounded by $\Phi_{11}$. A known norm bound on each individual disturbance $w_j$, $\|w_j\| \le \varepsilon$, implies a bound in the form of \eqref{equ:6} with $\Phi_{11} = (J \varepsilon^2) I$, $\Phi_{12} = 0$, and $\Phi_{22} = -I$. 

Plugging $W = X^+ - A X - B U$ into \eqref{equ:6} and after some algebra, we obtain 
\small \vspace{-6mm}
\begin{align}\label{equ:7}
\begin{bmatrix}
I \\ A^{\top} \\ B^{\top} 
\end{bmatrix}^{\!\!\top} \!\! \begin{bmatrix}
\Theta & \star & \star \\ - X \Phi_{12}^{\top} - X \Phi_{22} (X^+)^{\top} & X \Phi_{22} X^{\top} & \star \\ - U \Phi_{12}^{\top} - U \Phi_{22} (X^+)^{\top} & U \Phi_{22} X^{\top} & U \Phi_{22} U^{\top}
\end{bmatrix} & \!\! \begin{bmatrix}
I \\ A^{\top} \\ B^{\top}
\end{bmatrix} \nonumber \\[1mm]
& \!\!\!\!\!\!\!\!\!\!\!\! \succeq 0
\end{align} \normalsize
where $\Theta = \Phi_{11} + X^+ \Phi_{12}^{\top} + \Phi_{12} (X^+)^{\top} + X^+ \Phi_{22} (X^+)^{\top}$, and $\star$ indicates the transpose of the related element below the diagonal. We let $\Sigma$ be the collection of $(A,B)$ satisfying~\eqref{equ:7}:
\begin{equation}\label{equ:8}
\Sigma = \{(A,B): \text{\eqref{equ:7} is satisfied}\}
\end{equation}
According to \eqref{equ:5}, $(A_{\text{o}},B_{\text{o}}) \in \Sigma$.

Under Assumption~1, we propose a control approach to achieving the three objectives below \eqref{equ:1} for the unknown system \eqref{equ:1}. Our approach is based on the following matrix S-lemma developed in \cite{van2020noisy}:

{\bf Lemma~1 \cite{van2020noisy}.} Let $M,N \in \mathbb{R}^{(2n+m) \times (2n+m)}$ be symmetric and partitioned as follows:
\begin{equation}\label{equ:9}
M = \begin{bmatrix} M_{11} & M_{12} \\ M_{12}^{\top} & M_{22} \end{bmatrix} \quad N = \begin{bmatrix} N_{11} & N_{12} \\ N_{12}^{\top} & N_{22} \end{bmatrix}
\end{equation}
where $M_{11},N_{11} \in \mathbb{R}^{n \times n}$. Suppose $M_{22} \preceq 0$, $N_{22} \preceq 0$, $\text{ker}(N_{22}) \subseteq \text{ker}(N_{12})$, and there exists $\bar{Z} \in \mathbb{R}^{(n+m) \times n}$ satisfying $\begin{bmatrix} I \\ \bar{Z} \end{bmatrix}^{\top} N \begin{bmatrix} I \\ \bar{Z} \end{bmatrix} \succ 0$. Then, we have that
\begin{equation}\label{equ:10}
\begin{bmatrix} I \\ Z \end{bmatrix}^{\top} \!M \begin{bmatrix} I \\ Z \end{bmatrix} \succ 0 \text{ for all $Z$ satisfying } \begin{bmatrix} I \\ Z \end{bmatrix}^{\top} \!N \begin{bmatrix} I \\ Z \end{bmatrix} \succeq 0
\end{equation}
if and only if there exist $\alpha \ge 0$ and $\beta > 0$ such that
\begin{equation}\label{equ:11}
M - \alpha N \succeq \begin{bmatrix} \beta I & 0 \\ 0 & 0 \end{bmatrix}
\end{equation}

{\bf Proof.} See Theorem~13 of \cite{van2020noisy}. $\blacksquare$ 

\section{Constrained $\mathcal{H}_{\infty}$ Control}\label{sec:3}

In this section, we consider a linear-feedback $u(t) = K x(t)$ with a constant gain matrix $K \in \mathbb{R}^{m \times n}$. This controller yields the following closed-loop system:
\begin{subequations}\label{equ:12}
\begin{align}
x(t+1) &= \left(A_{\text{o}} + B_{\text{o}} K\right) x(t) + w(t) \label{equ:12_1} \\
y_1(t) &= \left(C_1 + D_1 K\right) x(t) \\
y_2(t) &= \left(C_2 + D_2 K\right) x(t)
\end{align}
\end{subequations}
The closed-loop transfer matrix from disturbance input $w$ to performance output $y_1$ is given by
\begin{equation}\label{equ:13}
G_1(z) = \left(C_1 + D_1 K\right) \left(zI - (A_{\text{o}} + B_{\text{o}} K)\right)^{-1}
\end{equation}
For a performance level $\gamma > 0$, the closed-loop matrix $A_{\text{c}} = A_{\text{o}} + B_{\text{o}} K$ is Schur stable and the $\mathcal{H}_{\infty}$ norm of $G_1(z)$ satisfies $\|G_1(z)\|_{\mathcal{H}_{\infty}} < \gamma$ if and only if there exists a matrix $P = P^{\top} \succ 0$ satisfying 
\small \vspace{-6mm}
\begin{subequations}\label{equ:14}
\begin{align}
& P - (A_{\text{o}} + B_{\text{o}} K)^{\top} P (A_{\text{o}} + B_{\text{o}} K) - (C_1 + D_1 K)^{\top} (C_1 + D_1 K) \nonumber \\
& - (A_{\text{o}} + B_{\text{o}} K)^{\top} P (\gamma^2 I - P)^{-1} P (A_{\text{o}} + B_{\text{o}} K) \succ 0  \\
& \gamma^2 I - P \succ 0
\end{align}
\end{subequations} \normalsize
See Theorem~2.2 of \cite{de1992discrete}. Now let $Q = P^{-1}$ (hence, $Q = Q^{\top} \succ 0$) and $Y = KQ$. With some algebra and Schur complement arguments, it can be shown that \eqref{equ:14} is equivalent to 
\small \vspace{-6mm}
\begin{subequations}\label{equ:15}
\begin{align}
& R = Q - (C_1 Q + D_1 Y)^{\top}(C_1 Q + D_1 Y) \succ 0 \label{equ:15_1} \\
& Q - \gamma^{-2} I - (A_{\text{o}} Q + B_{\text{o}} Y) R^{-1}(A_{\text{o}} Q + B_{\text{o}} Y)^{\top} \succ 0 \label{equ:15_2}
\end{align}
\end{subequations} \normalsize
Note that \eqref{equ:15_2} can be written as
\begin{equation}\label{equ:16}
\begin{bmatrix}
I \\ A_{\text{o}}^{\top} \\ B_{\text{o}}^{\top}
\end{bmatrix}^{\!\top}  \!\! \begin{bmatrix} Q - \gamma^{-2} I & 0 & 0 \\ 
0 & - Q R^{-1} Q & - Q R^{-1} Y^{\top} \\
0 & - Y R^{-1} Q & - Y R^{-1} Y^{\top} 
\end{bmatrix} \begin{bmatrix}
I \\ A_{\text{o}}^{\top} \\ B_{\text{o}}^{\top}
\end{bmatrix} \succ 0
\end{equation}
Now consider the quadratic Lyapunov function $V(x(t)) = x(t)^{\top} P x(t)$. When \eqref{equ:14} holds, we can derive the following dissipation inequality: 
\small \vspace{-6mm}
\begin{align}\label{equ:17}
& V(x(t+1)) - V(x(t)) = x(t+1)^{\top} P x(t+1) - x(t)^{\top} P x(t) \nonumber \\ 
&= x(t)^{\top} (A_{\text{c}}^{\top} P A_{\text{c}} - P) x(t) + w(t)^{\top} P w(t) \nonumber \\
&\,\, + x(t)^{\top} A_{\text{c}}^{\top} P w(t) + w(t)^{\top} P A_{\text{c}} x(t) \nonumber \\
&\le - x(t)^{\top} (C_1 + D_1 K)^{\top} (C_1 + D_1 K) x(t) + \gamma^2 w(t)^{\top} w(t) \nonumber \\
&\,\, - x(t)^{\top} A_{\text{c}}^{\top} P (\gamma^2 I - P)^{-1} P A_{\text{c}} x(t) - w(t)^{\top} (\gamma^2 I - P) w(t) \nonumber \\
&\,\, + x(t)^{\top} A_{\text{c}}^{\top} P w(t) + w(t)^{\top} P A_{\text{c}} x(t) \nonumber \\
&= - y_1(t)^{\top} y_1(t) + \gamma^2 w(t)^{\top} w(t) - S \nonumber \\
&\le - \|y_1(t)\|^2 + \gamma^2 \|w(t)\|^2
\end{align} \normalsize
where $S = \|(\gamma^2 I - P)^{-1/2} P A_{\text{c}} x(t) - (\gamma^2 I - P)^{1/2} w(t)\|^2$ satisfies $S \ge 0$. We define, for $r > 0$, the ellipsoidal set
\begin{equation}\label{equ:18}
\mathcal{E}(P,r) = \{x \in \mathbb{R}^n: V(x) \le r\}
\end{equation}
and we arrive at the following result:

{\bf Lemma~2.} Suppose \eqref{equ:14} holds and the energy of the disturbance input is bounded as $\sum_{t = 0}^{\infty} \|w(t)\|^2 \le \sigma_0$ for some $\sigma_0 \ge 0$. Then, for any $r_0 \ge x(0)^{\top} P x(0) + \gamma^2 \sigma_0$, $\mathcal{E}(P,r_0)$ is an invariant set of \eqref{equ:12}, i.e., $x(t) \in \mathcal{E}(P,r_0)$ for all $t \in \mathbb{N}$. 

{\bf Proof.} Suppose \eqref{equ:14} holds, $\sum_{t = 0}^{\infty} \|w(t)\|^2 \le \sigma_0$, and $r_0 \ge x(0)^{\top} P x(0) + \gamma^2 \sigma_0$. Then, the dissipation inequality \eqref{equ:17} implies
\begin{align}
& V(x(t)) \le V(x(0)) - \sum_{i=0}^{t-1} \|y_1(i)\|^2 + \gamma^2 \sum_{i=0}^{t-1} \|w(i)\|^2 \nonumber \\
&\quad\quad \le x(0)^{\top} P x(0) + \gamma^2 \sum_{i=0}^{\infty} \|w(i)\|^2 \le r_0
\end{align}
for all $t \in \mathbb{N}$. This shows $x(t) \in \mathcal{E}(P,r_0)$ for all $t \in \mathbb{N}$. $\blacksquare$

We now consider the optimization problem \eqref{equ:20} for designing the feedback gain $K = Y Q^{-1}$. In \eqref{equ:20}, $e_v$ denotes the $v$th standard basis vector, $\sigma_0 > 0$ is a forecasted energy bound of the disturbance input, and $r_0 > 0$ is a design parameter to be tuned so that \eqref{equ:20} is feasible. We note that \eqref{equ:20} is an LMI (hence, convex) problem in the decision variables $(\eta, Q, Y, \alpha, \beta)$. A design of $K$ based on~\eqref{equ:20} has several desirable properties, stated in the following theorem:

\begin{table*}
\centering
\begin{minipage}{1\textwidth}
\hrule
\begin{subequations}\label{equ:20}
\begin{align}
& \max_{\eta > 0, Q = Q^{\top}\!, Y, \alpha \ge 0, \beta > 0} \, \eta \label{equ:20_1} \\
\text{s.t. } & \begin{bmatrix} Q - (\eta+\beta) I & \star & \star & \star & \star \\ 
0 & 0 & \star & \star & \star \\
0 & 0 & 0 & \star & \star \\
0 & Q & Y^{\top} & Q & \star \\
0 & 0 & 0 & C_1 Q + D_1 Y & I \end{bmatrix} \succeq \alpha \begin{bmatrix}
\Theta & \star & \star & \star & \,\star \\ 
- X \Phi_{12}^{\top} - X \Phi_{22} (X^+)^{\top} & X \Phi_{22} X^{\top} & \star & \star & \,\star \\ 
- U \Phi_{12}^{\top} - U \Phi_{22} (X^+)^{\top} & U \Phi_{22} X^{\top} & U \Phi_{22} U^{\top} & \star & \,\star \\
0 & 0 & 0 & 0 & \,\star \\
0 & 0 & 0 & 0 & \,0
\end{bmatrix} \label{equ:20_2} \\
& \begin{bmatrix} Q & \star \\ C_1 Q + D_1 Y & I \end{bmatrix} \succ 0 \quad\,(\stepcounter{equation}{\text{\theequation}})\quad\,\,\,
\begin{bmatrix} r_0 & \star & \star \\ x(0) & Q & \star \\ 1 & 0 & \sigma_0^{-1} \eta \end{bmatrix} \succeq 0 \quad\,(\stepcounter{equation}{\text{\theequation}})\quad\,\,\,
\begin{bmatrix} y_{2v,\max}^2\, r_0^{-1} & \star \\ (C_2 Q + D_2 Y)^{\top}e_v & Q \end{bmatrix} \succeq 0,\,\, v = 1,2,\dots,p_2 \label{equ:20_5}
\end{align}
\end{subequations}
where $\Theta = \Phi_{11} + X^+ \Phi_{12}^{\top} + \Phi_{12} (X^+)^{\top} + X^+ \Phi_{22} (X^+)^{\top}$, and $\star$ indicates the transpose of the related element below the diagonal.
\medskip
\hrule
\end{minipage}
\end{table*}

{\bf Theorem~1.} Suppose 
\begin{enumerate}
    \item[(a)] The offline data satisfy Assumption~1 and there exists $(A,B) = (\bar{A},\bar{B})$ such that \eqref{equ:7} is strictly positive-definite;
    \item[(b)] The online disturbance inputs satisfy the energy bound $\sum_{t = 0}^{\infty} \|w(t)\|^2 \le \sigma_0$;
    \item[(c)] The tuple $(\eta_0, Q_0, Y_0, \alpha_0, \beta_0)$ is a solution to~\eqref{equ:20}.
\end{enumerate}
Then, if we control the unknown system \eqref{equ:1} using the linear-feedback $u(t) = K_0 x(t)$, with $K_0 = Y_0 Q_0^{-1}$, the closed-loop system has the following properties:
\begin{enumerate}
    \item[(i)] The closed-loop matrix $A_{\text{c},0} = A_{\text{o}} + B_{\text{o}} K_0$ is Schur stable and the closed-loop $\mathcal{H}_{\infty}$ norm from $w$ to $y_1$ is less than $\gamma_0 = \eta_0^{-1/2}$;
    \item[(ii)] The constraints in \eqref{equ:2} are satisfied at all times.
\end{enumerate}

{\bf Proof.} First, let \small
\begin{align}\label{equ:21}
& M = \left[\begin{array}{c|c} M_{11} & M_{12} \\ \hline M_{12}^{\top} & M_{22} \end{array}\right] = \left[\begin{array}{c|cc}
Q - \eta I & 0 & 0 \\ \hline
0 & - Q R^{-1} Q & - Q R^{-1} Y^{\top} \\
0 & - Y R^{-1} Q & - Y R^{-1} Y^{\top} 
\end{array}\right] \nonumber \\[2mm]
& N = \left[\begin{array}{c|c} N_{11} & N_{12} \\ \hline N_{12}^{\top} & N_{22} \end{array}\right] \nonumber \\
&= \left[\begin{array}{c|cc}
\Theta & \star & \star \\ \hline
- X \Phi_{12}^{\top} - X \Phi_{22} (X^+)^{\top} & X \Phi_{22} X^{\top} & \star \\ - U \Phi_{12}^{\top} - U \Phi_{22} (X^+)^{\top} & U \Phi_{22} X^{\top} & U \Phi_{22} U^{\top}
\end{array}\right] 
\end{align} \normalsize
where $R = Q - (C_1 Q + D_1 Y)^{\top}(C_1 Q + D_1 Y)$. Using a Schur complement argument, it can be shown that \eqref{equ:20_2} is equivalent to
\begin{equation}\label{equ:22}
M - \alpha N \succeq \begin{bmatrix} \beta I & 0 \\ 0 & 0 \end{bmatrix}
\end{equation}
Meanwhile, for the $M$ and $N$ defined and partitioned as above, the following conditions can be verified: 1) $M_{22} \preceq 0$, using the fact that $R \succ 0$ due to the constraint (20c), 2)~$N_{22} \preceq 0$, using $\Phi_{22} \prec 0$ according to Assumption~1, and 3) $\text{ker}(N_{22}) \subseteq \text{ker}(N_{12})$. Now with $Z = [A, B]^{\top}$ and $\bar{Z} = [\bar{A}, \bar{B}]^{\top}$ (with $(\bar{A},\bar{B})$ given in assumption~(a)), it can be seen that the assumptions of Lemma~1 are all satisfied. According to Lemma~1, \eqref{equ:22} holds for some $\alpha \ge 0$ and $\beta > 0$ if and only if
\begin{equation}\label{equ:23}
\begin{bmatrix} I \\ A^{\top} \\ B^{\top} \end{bmatrix}^{\top} M \begin{bmatrix} I \\ A^{\top} \\ B^{\top} \end{bmatrix} \succ 0 \text{ for all $(A,B) \in \Sigma$}
\end{equation}
where $\Sigma$ is defined in \eqref{equ:8}. Recall that $(A_{\text{o}},B_{\text{o}}) \in \Sigma$. Therefore, we have shown that if $(\eta, Q = Q^{\top}, Y, \alpha, \beta)$ satisfies \eqref{equ:20_2}, (20c), and $\eta > 0, \alpha \ge 0, \beta > 0$, then~\eqref{equ:16} (hence, \eqref{equ:15_2}), with $\gamma = \eta^{-1/2}$, necessarily holds. Note that (20c) also implies \eqref{equ:15_1} and $Q \succ 0$. Now, if we let $P = Q^{-1}$ and $K = Y Q^{-1}$, we arrive at \eqref{equ:14} with $\gamma = \eta^{-1/2}$. This proves part~(i).

Now suppose $(\eta, Q = Q^{\top}, Y, \alpha, \beta)$ satisfies \eqref{equ:20_2}-\eqref{equ:20_5} and $\eta > 0, \alpha \ge 0, \beta > 0$. Above we have shown that, with $P = Q^{-1}$, $K = Y Q^{-1}$, and $\gamma = \eta^{-1/2}$, \eqref{equ:14} holds. Meanwhile, using a Schur complement argument, (20d) implies $r_0 \ge x(0)^{\top} Q^{-1} x(0) + \eta^{-1} \sigma_0 = x(0)^{\top} P x(0) + \gamma^2 \sigma_0$. In this case, according to Lemma~2, $\sum_{t = 0}^{\infty} \|w(t)\|^2 \le \sigma_0$ in assumption~(b) implies that $\mathcal{E}(P,r_0)$ is an invariant set of the closed-loop system, i.e., $x(t) \in \mathcal{E}(P,r_0)$ for all $t \in \mathbb{N}$. Recall that $\mathcal{E}(P,r)$ is the ellipsoidal set defined in \eqref{equ:18}. The support function of $\mathcal{E}(P,r)$ is $h_{\mathcal{E}(P,r)}(\zeta) = \sqrt{r\, \zeta^{\top} P^{-1} \zeta}$. Hence, $x(t) \in \mathcal{E}(P,r_0)$ implies
\begin{align}\label{equ:24}
&\! y_{2v}(t) = e_v^{\top}(C_2 + D_2 K) x(t) \le h_{\mathcal{E}(P,r_0)}((C_2 + D_2 K)^{\top}e_v) \nonumber \\
&= \sqrt{r_0\, e_v^{\top}(C_2 + D_2 K) P^{-1} (C_2 + D_2 K)^{\top}e_v} \nonumber \\
&= \sqrt{r_0\, e_v^{\top}(C_2 Q + D_2 Y) Q^{-1} (C_2 Q + D_2 Y)^{\top}e_v}
\end{align}
Meanwhile, using a Schur complement argument, \eqref{equ:20_5} is equivalent to 
\small \vspace{-6mm}
\begin{align}\label{equ:25}
& y_{2v,\max}^2\, r_0^{-1} - e_v^{\top}(C_2 Q + D_2 Y) Q^{-1} (C_2 Q + D_2 Y)^{\top}e_v \ge 0 \nonumber \\[1mm]
& \iff \nonumber \\[1mm]
& \sqrt{r_0\, e_v^{\top}(C_2 Q + D_2 Y) Q^{-1} (C_2 Q + D_2 Y)^{\top}e_v} \le y_{2v,\max}
\end{align} \normalsize
Combing \eqref{equ:24} and \eqref{equ:25}, we obtain $y_{2v}(t) \le y_{2v,\max}$. This completes the proof of part~(ii). $\blacksquare$

Theorem~1 states the stability, disturbance attenuation, and constraint enforcement properties of the feedback gain $K$ for the unknown system \eqref{equ:1} synthesized using its trajectory data according to \eqref{equ:20}. The existence of $(\bar{A},\bar{B})$ such that \eqref{equ:7} is strictly positive-definite in assumption~(a) of Theorem~1 can be checked offline. We make the following two remarks about Theorem~1:

{\bf Remark~1.} The $\mathcal{H}_{\infty}$ norm from $w$ to $y_1$ represents a level of disturbance attenuation of the closed-loop system. In particular, using \eqref{equ:17} recursively we can obtain the following dissipation inequality that holds for all $t \in \mathbb{N}$:
\begin{equation}\label{equ:26}
\sum_{i=0}^t \|y_1(i)\|^2 \le \gamma^2 \sum_{i=0}^t \|w(i)\|^2 + x(0)^{\top} P x(0)
\end{equation}
and this indicates that the $\l_2$-gain from disturbance $w$ to output $y_1$ is bounded by $\gamma$. Because \eqref{equ:20} maximizes $\eta > 0$, which is equivalent to minimizing $\gamma = \eta^{-1/2}$, the obtained controller seeks to maximize disturbance attenuation while satisfying the time-domain constraints~in~\eqref{equ:2}.

{\bf Remark~2.} Differently from many other robust control formulations that assume set-bounded disturbances (i.e., $w(t) \in \mathbb{W}$ for some known set bound $\mathbb{W}$ and for all $t \in \mathbb{N}$), our formulation \eqref{equ:20} enforces constraints for disturbance inputs that have bounded total energy where the energy can be distributed arbitrarily over time (i.e., $\sum_{t = 0}^{\infty} \|w(t)\|^2 \le \sigma_0$). This total energy model is particularly suitable for modeling transient disturbances such as wind gusts in aircraft flight control or wind turbine control, power outages in power systems, temporary actuator or sensor failures, etc. Such disturbances occur infrequently and typically last a short period of time (as compared to persistent disturbances), but they can have a significant magnitude. Meanwhile, predicting exactly when they will occur can be difficult or impossible. In such a case, on the one hand, our formulation based on a total energy disturbance model may lead to a less conservative solution (hence, better performance) than those assuming a set bound at all times; on the other hand, estimating/predicting an energy bound for such transient disturbance events using historic data and/or real-time information (such as weather conditions) is possible in many applications. In what follows we assume a mechanism that is able to forecast an energy bound for future disturbances is available. Treating set-bounded persistent disturbances is left for future research.

\section{Moving-Horizon Control}\label{sec:4}

We now consider a strategy for implementing the control developed in Section~\ref{sec:3} in a moving-horizon manner: At~each time $t \in \mathbb{N}$, one solves \eqref{equ:20} with the current system state as the initial condition $x(0)$ in (20d) for a feedback gain $K_t$ and uses the control $u(t) = K_t x(t)$ over one time step. This way, the feedback gain becomes adaptive to system state and the control becomes nonlinear, possibly leading to improved performance while satisfying constraints. However, this simple implementation may fail to guarantee stability and disturbance attenuation, as discussed in \cite{chen2004disturbance,chen2006moving}. To recover a disturbance attenuation guarantee, following the strategy of \cite{chen2004disturbance,chen2006moving}, we consider the following inequality:
\begin{equation}\label{equ:27}
x(t)^{\top} P_t x(t) - x(t)^{\top} P_{t-1} x(t) \le \Delta_t
\end{equation}
where $P_t$ is associated with the feedback gain $K_t$ that is used at time $t$ and with an $\mathcal{H}_{\infty}$ performance level of $\gamma_t$ (i.e., the triple $(P_t, K_t, \gamma_t)$ satisfies \eqref{equ:14}), $P_{t-1}$ is associated with $K_{t-1}$ and $\gamma_{t-1}$, and $\Delta_t$ keeps track of a previous dissipation level and is defined recursively according to
\begin{equation}\label{equ:28}
\! \Delta_t \!=\! \Delta_{t-1} - \Big(x(t-1)^{\top} P_{t-1} x(t-1) - x(t-1)^{\top} P_{t-2} x(t-1)\Big)
\end{equation}
for $t \ge 2$ and $\Delta_1 = 0$. Note that the definition of $\Delta_t$ in \eqref{equ:28} only uses state and $P$-matrix information up to time $t-1$.

{\bf Lemma~3.} Suppose \eqref{equ:27} holds for all $t \ge 1$ where $\Delta_t$ is defined recursively according to \eqref{equ:28}. Then, the following dissipation inequality will be satisfied for all $t \in \mathbb{N}$:
\begin{equation}\label{equ:29}
\sum_{i=0}^t \|y_1(i)\|^2 \le \bar{\gamma}_t^2 \sum_{i=0}^t \|w(i)\|^2 + x(0)^{\top} P_0 x(0)
\end{equation}
where $\bar{\gamma}_t = \max\{\gamma_0,\gamma_1,\dots,\gamma_t\}$.

{\bf Proof.} For each $i$, because $(P_i, K_i, \gamma_i)$ satisfies \eqref{equ:14}, similar to \eqref{equ:17}, we can derive the following inequality:
\begin{equation}\label{equ:30}
x(i+1)^{\top} P_i x(i+1) - x(i)^{\top} P_i x(i) \le - \|y_1(i)\|^2 + \gamma_i^2 \|w(i)\|^2
\end{equation}
Summing over $i = 0,1,\dots,t$, we obtain
\begin{align}\label{equ:31}
&\! \sum_{i = 0}^t \|y_1(i)\|^2 \le \sum_{i = 0}^t \gamma_i^2 \|w(i)\|^2 + x(0)^{\top} P_0 x(0) \,+\, \\
&\! \sum_{i = 1}^t \!\Big(x(i)^{\top} P_i x(i) -  x(i)^{\top} P_{i-1} x(i)\Big) \!- x(t+1)^{\top} P_t x(t+1) \nonumber
\end{align}
The definition \eqref{equ:28} yields the following closed-form expression for $\Delta_t$:
\begin{equation}\label{equ:32}
\Delta_t = - \sum_{i = 1}^{t-1} \Big(x(i)^{\top} P_i x(i) - x(i)^{\top} P_{i-1} x(i)\Big)
\end{equation}
Supposing \eqref{equ:27} holds at $t$, we have
\begin{equation}\label{equ:33}
\sum_{i = 1}^t \Big(x(i)^{\top} P_i x(i) - x(i)^{\top} P_{i-1} x(i)\Big) \le 0
\end{equation}
Combining \eqref{equ:31}, \eqref{equ:33}, and $\bar{\gamma}_t = \max_{i=0,\dots,t} \gamma_i$, we obtain
\begin{align}\label{equ:34}
\sum_{i = 0}^t \|y_1(i)\|^2 &\le \sum_{i = 0}^t \gamma_i^2 \|w(i)\|^2 + x(0)^{\top} P_0 x(0) \nonumber \\
&\quad\quad\quad - x(t+1)^{\top} P_t x(t+1) \nonumber \\
&\le \bar{\gamma}_t^2 \sum_{i = 0}^t \|w(i)\|^2 + x(0)^{\top} P_0 x(0) \quad \blacksquare
\end{align}

We now consider the following moving-horizon approach to constrained $\mathcal{H}_{\infty}$ control for the unknown system \eqref{equ:1}: At~each time $t \in \mathbb{N}$, we solve the optimization problem~\eqref{equ:35} for designing the feedback gain $K_t = Y_t Q_t^{-1}$ and use the control $u(t) = K_t x(t)$ over one time step. In particular, the constraint \eqref{equ:35_6} is excluded at the initial time $t = 0$ and included for $t \ge 1$. In \eqref{equ:35}, $x(t)$ denotes the measured current system state, $P_{t-1} = Q_{t-1}^{-1}$ is from the previous time, $\Delta_t$ is defined according to \eqref{equ:28}, $\sigma_t$ represents a forecasted bound on the total energy of present and future disturbances, i.e., $\sum_{i = t}^{\infty} \|w(i)\|^2 \le \sigma_t$, and $r_t > 0$ is a design parameter, of which a design method is informed by Lemma~4 and Theorem~2. 
For a moving-horizon control algorithm, recursive feasibility (i.e., the online optimization problem being feasible at a given time implies the problem being feasible again at the next time) is a highly desirable property. Before we discuss the closed-loop properties of the proposed algorithm, the following lemma provides a recursive feasibility result:

\begin{table*}
\centering
\begin{minipage}{1\textwidth}
\hrule
\begin{subequations}\label{equ:35}
\begin{align}
& \max_{\eta > 0, Q = Q^{\top}\!, Y, \alpha \ge 0, \beta > 0} \, \eta \label{equ:35_1} \\
\text{s.t. } & \begin{bmatrix} Q - (\eta+\beta) I & \star & \star & \star & \star \\ 
0 & 0 & \star & \star & \star \\
0 & 0 & 0 & \star & \star \\
0 & Q & Y^{\top} & Q & \star \\
0 & 0 & 0 & C_1 Q + D_1 Y & I \end{bmatrix} \succeq \alpha \begin{bmatrix}
\Theta & \star & \star & \star & \,\star \\ 
- X \Phi_{12}^{\top} - X \Phi_{22} (X^+)^{\top} & X \Phi_{22} X^{\top} & \star & \star & \,\star \\ 
- U \Phi_{12}^{\top} - U \Phi_{22} (X^+)^{\top} & U \Phi_{22} X^{\top} & U \Phi_{22} U^{\top} & \star & \,\star \\
0 & 0 & 0 & 0 & \,\star \\
0 & 0 & 0 & 0 & \,0
\end{bmatrix} \label{equ:35_2} \\
& \begin{bmatrix} Q & \star \\ C_1 Q + D_1 Y & I \end{bmatrix} \succ 0 \quad\quad\quad\quad\quad\quad\quad\quad\,(\stepcounter{equation}{\text{\theequation}})\quad\quad\quad\quad\quad\quad\quad \begin{bmatrix} r_t & \star & \star \\ x(t) & Q & \star \\ 1 & 0 & \sigma_t^{-1} \eta \end{bmatrix} \succeq 0
\label{equ:35_4} \\
& \begin{bmatrix} y_{2v,\max}^2\, r_t^{-1} & \star \\ (C_2 Q + D_2 Y)^{\top}e_v & Q \end{bmatrix} \succeq 0,\,\, v = 1,2,\dots,p_2 \quad\quad (\stepcounter{equation}{\text{\theequation}})\quad\quad\quad\,
\begin{bmatrix} x(t)^{\top} P_{t-1} x(t) + \Delta_t & \star \\ x(t) & Q \end{bmatrix} \succeq 0
\label{equ:35_6}
\end{align}
\end{subequations}
where $\Theta = \Phi_{11} + X^+ \Phi_{12}^{\top} + \Phi_{12} (X^+)^{\top} + X^+ \Phi_{22} (X^+)^{\top}$, $\star$ indicates the transpose of the related element below the diagonal, and the constraint \eqref{equ:35_6} is excluded at $t = 0$ and included for $t \ge 1$.
\medskip
\hrule
\end{minipage}
\end{table*}

{\bf Lemma~4.} Suppose \eqref{equ:35} is feasible at time $t$ and $(\eta_t, Q_t, Y_t, \alpha_t, \beta_t)$ denotes a solution. At time $t+1$, if~the forecasted disturbance energy bound $\sigma_{t+1}$ satisfies $\sigma_{t+1} \le \sigma_t - \|w(t)\|^2$ and the parameter $r_{t+1}$ is chosen to be $r_{t+1} = r_t$, \eqref{equ:35} will be feasible again. In particular, in this case, $(\eta_t, Q_t, Y_t, \alpha_t, \beta_t)$ remains to be a solution to \eqref{equ:35} at time $t+1$.

{\bf Proof.} First, we note that the constraints \eqref{equ:35_2}, (35c), and (35e) with $r_{t+1} = r_t$ do not change from $t$ to $t+1$. The solution $(\eta_t, Q_t, Y_t, \alpha_t, \beta_t)$ satisfying \eqref{equ:35_2} and (35c) implies \eqref{equ:14} with $P = Q_t^{-1}$, $K = Y_t Q_t^{-1}$, and $\gamma = \eta_t^{-1/2}$. In this case, similar to \eqref{equ:17}, the following inequality holds
\begin{align}\label{equ:36}
& x(t+1)^{\top} Q_t^{-1} x(t+1) - x(t)^{\top} Q_t^{-1} x(t) \nonumber \\
&\quad \le - \|y_1(t)\|^2 + \eta_t^{-1} \|w(t)\|^2
\end{align}
Meanwhile, $(\eta_t, Q_t)$ satisfying the constraint \eqref{equ:35_4} at time $t$ is equivalent to
\begin{equation}\label{equ:37}
r_t \ge x(t)^{\top} Q_t^{-1} x(t) + \eta_t^{-1} \sigma_t    
\end{equation}
If $\sigma_{t+1} \le \sigma_t - \|w(t)\|^2$ and $r_{t+1} = r_t$, together with \eqref{equ:36} and \eqref{equ:37}, we can obtain
\begin{align}\label{equ:38}
& x(t+1)^{\top} Q_t^{-1} x(t+1) + \eta_t^{-1} \sigma_{t+1} \nonumber \\
&\le x(t)^{\top} Q_t^{-1} x(t) + \eta_t^{-1} \|w(t)\|^2 + \eta_t^{-1} (\sigma_t - \|w(t)\|^2) \nonumber \\
&= x(t)^{\top} Q_t^{-1} x(t) + \eta_t^{-1} \sigma_t \nonumber \\
&\le r_t = r_{t+1}
\end{align}
which implies that $(\eta_t, Q_t)$ also satisfies the constraint~\eqref{equ:35_4} at time $t+1$. Last, for $t = 0$, we have $\Delta_{t+1} = \Delta_1 = 0$; for $t \ge 1$, $Q_t$ satisfying the constraint~\eqref{equ:35_6} implies \eqref{equ:27} with $P_t = Q_t^{-1}$. In the latter case, similar to \eqref{equ:32}--\eqref{equ:33}, we have
\begin{equation}\label{equ:39}
\Delta_{t+1} = -\sum_{i = 1}^t \Big(x(i)^{\top} P_i x(i) - x(i)^{\top} P_{i-1} x(i)\Big) \ge 0
\end{equation}
At time $t+1$, the constraint \eqref{equ:35_6} is equivalent to
\begin{equation}\label{equ:40}
x(t+1)^{\top} Q^{-1} x(t+1) - x(t+1)^{\top} P_t x(t+1) \le \Delta_{t+1}
\end{equation}
where $P_t = Q_t^{-1}$. Because $\Delta_{t+1} \ge 0$, it is clear that $Q = Q_t$ satisfies \eqref{equ:40} (hence, \eqref{equ:35_6}).

Therefore, we have shown that if $\sigma_{t+1} \le \sigma_t - \|w(t)\|^2$ and $r_{t+1} = r_t$, a solution at time $t$, $(\eta_t, Q_t, Y_t, \alpha_t, \beta_t)$, still satisfies all constraints of \eqref{equ:35} at time $t+1$, i.e., $(\eta_t, Q_t, Y_t, \alpha_t, \beta_t)$ remains to be a solution to \eqref{equ:35} at $t+1$. This proves the result. $\blacksquare$

{\bf Remark~3.} Lemma~4 represents a sufficient condition for the online optimization problem \eqref{equ:35} to be recursively feasible. The assumption $\sigma_{t+1} \le \sigma_t - \|w(t)\|^2$ is reasonable because $\sigma_t$ bounds the total energy of disturbance inputs from time $t$ and $\sigma_{t+1}$ bounds that from $t + 1$ -- they differ by the energy of the disturbance input at time~$t$. Then, $r_{t+1} = r_t$ yields a simple strategy for setting the parameter $r_t$ -- at each time $t$, $r_t$ is set to its previous value. At time instants at which \eqref{equ:35} is infeasible with $r_t = r_{t-1}$ (due to errors in forecasted disturbance energy bound and violations of $\sigma_{t+1} \le \sigma_t - \|w(t)\|^2$), an alternative strategy that promotes feasibility is to include $r_t$ or its inverse $\rho_t = r_t^{-1}$ as a decision variable optimized together with the other variables $(\eta, Q, Y, \alpha, \beta)$. In~this case, \eqref{equ:35} becomes a nonconvex problem with a single nonconvex variable $r_t$ or $\rho_t$, which can be solved using a branch-and-bound type algorithm with branching on~$r_t$~or~$\rho_t$~\cite{tuy1992nonconvex}.

We now analyze the closed-loop properties of the system under the proposed moving-horizon control approach. The properties are given in Theorem~2.

{\bf Theorem~2.} Suppose 
\begin{enumerate}
    \item[(a)] The offline data satisfy Assumption~1 and there exists $(A,B) = (\bar{A},\bar{B})$ such that \eqref{equ:7} is strictly positive-definite;
    \item[(b)] The online disturbance inputs satisfy $\sum_{i = t}^{\infty} \|w(i)\|^2 \le \sigma_t$ for all $t \in \mathbb{N}$, where $\sigma_t$ are forecasted energy bounds and used in \eqref{equ:35};
    \item[(c)] The online optimization problem \eqref{equ:35} is feasible at all $t \in \mathbb{N}$ and $(\eta_t, Q_t, Y_t, \alpha_t, \beta_t)$ denotes a solution to \eqref{equ:35} at each $t$;
    \item[(d)] The solutions $(\eta_t, Q_t, Y_t, \alpha_t, \beta_t)$ have a common lower bound $\underline{\eta} > 0$ on their objective values, i.e., $\eta_t \ge \underline{\eta}$ for all $t \in \mathbb{N}$.
\end{enumerate}
Then, if we control the unknown system \eqref{equ:1} using $u(t) = K_t x(t)$, with $K_t = Y_t Q_t^{-1}$, at all $t \in \mathbb{N}$, the closed-loop system has the following properties:
\begin{enumerate}
    \item[(i)] The system state $x(t)$ converges to $0$ as $t \to \infty$;
    \item[(ii)] The following dissipation inequality is satisfied for all $t \in \mathbb{N}$:
    \begin{equation}\label{equ:41}
    \sum_{i=0}^t \|y_1(i)\|^2 \le \bar{\gamma}^2 \sum_{i=0}^t \|w(i)\|^2 + x(0)^{\top} P_0 x(0)
    \end{equation}
    where $\bar{\gamma} = \underline{\eta}^{-1/2}$, indicating that the $\l_2$-gain from disturbance $w$ to output $y_1$ is bounded by $\bar{\gamma}$;
    \item[(iii)] The constraints in \eqref{equ:2} are satisfied at all times.
\end{enumerate}
Furthermore, suppose (a) and (b) hold and
\begin{enumerate}
    \item[(e)] Problem \eqref{equ:35} is feasible at the initial time $t = 0$, the energy bounds $\sigma_t$ satisfy $\sigma_{t+1} \le \sigma_t - \|w(t)\|^2$ for all $t \in \mathbb{N}$, and the parameters $r_t$ are chosen as $r_t = r_0$ for all $t \in \mathbb{N}$;
    \item[(f)] The tuple $(\eta_t, Q_t, Y_t, \alpha_t, \beta_t)$ is a (not only feasible but also) optimal solution to \eqref{equ:35} at each $t \in \mathbb{N}$.
\end{enumerate}
Then, the following results hold:
\begin{enumerate}
    \item[(iv)] The conditions (c) and (d) are necessarily satisfied;
    \item[(v)] The solutions have non-decreasing objective values, i.e., $\eta_{t+1} \ge \eta_t$ for all $t \in \mathbb{N}$;
    \item[(vi)] The dissipation inequality \eqref{equ:41} holds true with the gain of $\bar{\gamma} = \eta_0^{-1/2}$.
\end{enumerate}

{\bf Proof.} We start with proving the inequality \eqref{equ:41} in (ii). At each time $t \in \mathbb{N}$, the solution $(\eta_t, Q_t, Y_t, \alpha_t, \beta_t)$ satisfies the constraints \eqref{equ:35_2} and (35c). Following similar steps as in the proof of Theorem~1, part~(i), we can show that \eqref{equ:14}, with $P = Q_t^{-1}$, $K = Y_t Q_t^{-1}$, and $\gamma = \eta_t^{-1/2}$, holds. The constraint \eqref{equ:35_6} is equivalent to
\begin{equation}\label{equ:42}
x(t)^{\top} Q^{-1} x(t) - x(t)^{\top} P_{t-1} x(t) \le \Delta_t
\end{equation}
Because \eqref{equ:35_6} (hence, \eqref{equ:42}) is satisfied by $Q_t$ at each $t \ge 1$, we have \eqref{equ:27} for all $t \ge 1$. In this case, according to Lemma~3, \eqref{equ:29} holds for all $t \in \mathbb{N}$. Meanwhile, $\eta_t \ge \underline{\eta}$ for all $t \in \mathbb{N}$ (assumption (d)) implies that $\bar{\gamma}_t \le \bar{\gamma}$ for all~$t$, where $\bar{\gamma}_t = \max_{i=0,\dots,t} \gamma_i$, $\gamma_i = \eta_i^{-1/2}$, and $\bar{\gamma} = \underline{\eta}^{-1/2}$. The combination of \eqref{equ:29} and $\bar{\gamma}_t \le \bar{\gamma}$ leads to \eqref{equ:41}. This proves part~(ii). Now, because \eqref{equ:41} holds for all $t \in \mathbb{N}$, as $t \to \infty$, we obtain
\begin{equation}\label{equ:43}
\sum_{i=0}^{\infty} \|y_1(i)\|^2 \le \bar{\gamma}^2 \sum_{i=0}^{\infty} \|w(i)\|^2 + x(0)^{\top} P_0 x(0)
\end{equation}
Note that under assumption (b), the right-hand side of~\eqref{equ:43} is bounded by $\bar{\gamma}^2 \sigma_0 + x(0)^{\top} P_0 x(0)$, and hence the series on the left-hand side converges according to the monotone convergence theorem. And \eqref{equ:43} implies $y_1(t) \to 0$ as $t \to \infty$. When $(C_1,A_{\text{o}})$ is detectable, $y_1(t) \to 0$ implies $x(t) \to 0$. This proves part~(i). For part~(iii), because $(\eta_t, Q_t, Y_t)$ satisfies \eqref{equ:35_4} and (35e), following similar steps as in the proof of Theorem~1, part~(ii), it can be shown that $x(t) \in \mathcal{E}(Q_t^{-1},r_t)$ (due to \eqref{equ:35_4}) and this implies $y_{2v}(t) = e_v^{\top}(C_2 + D_2 K_t) x(t) \le y_{2v,\max}$ for all $v = 1,2,\dots,p_2$ (due to (35e)).

Now suppose (a), (b), (e), and (f) hold. According to Lemma~4, when \eqref{equ:35} is feasible at $t = 0$, $\sigma_1 \le \sigma_0 - \|w(0)\|^2$, and $r_1 = r_0$, an optimal solution $(\eta_0, Q_0, Y_0, \alpha_0, \beta_0)$ to~\eqref{equ:35} at $t = 0$ remains to be a feasible solution to \eqref{equ:35} at $t = 1$. In this case, \eqref{equ:35} is feasible at $t = 1$ and an optimal solution $(\eta_1, Q_1, Y_1, \alpha_1, \beta_1)$ has an objective value at least as large as that of the feasible solution $(\eta_0, Q_0, Y_0, \alpha_0, \beta_0)$, i.e., $\eta_1 \ge \eta_0$. Using the same argument recursively, we can conclude that \eqref{equ:35} is feasible at all $t \in \mathbb{N}$ and $\eta_{t+1} \ge \eta_t$ for all $t$. The latter also implies that $\eta_0 > 0$ is a common lower bound for all $\eta_t$, i.e., (d) is satisfied with $\underline{\eta} = \eta_0$. Accordingly, \eqref{equ:41} holds with $\bar{\gamma} = \underline{\eta}^{-1/2} = \eta_0^{-1/2}$. This completes the proofs of parts~(iv), (v), and (vi). $\blacksquare$

We make the following remark about Theorem~2:

{\bf Remark~4.} Theorem~2 shows that our objectives of closed-loop stability, disturbance attenuation, and constraint enforcement stated in Section~\ref{sec:2} are fulfilled by the proposed data-driven moving-horizon $\mathcal{H}_{\infty}$ control based on~\eqref{equ:35}. In particular, part~(i) shows that the system state converges to zero asymptotically for any disturbance input signal that has bounded total energy. Parts~(v) and (vi) show that, under assumption~(e) (which is reasonable as discussed in Remark~3), the moving-horizon approach based on~\eqref{equ:35} achieves a disturbance attenuation performance at least as good as that achieved by a constant linear-feedback designed based on~\eqref{equ:20}. In particular, the moving-horizon approach based on~\eqref{equ:35} is able to dynamically adapt the feedback gain $K_t$ to measured/estimated system state $x(t)$ and forecasted disturbance energy bound $\sigma_t$ and hence has the potential to achieve improved disturbance attenuation performance while satisfying constraints. We will demonstrate this improvement with a numerical example in the following section.

\section{Numerical Example}\label{sec:5}

In this section, we use a numerical example to illustrate our proposed control approach. Consider a system in the form of \eqref{equ:1} with the following parameters:
\small \vspace{-6mm}
\begin{subequations}\label{equ:44} 
\begin{align}
A_{\text{o}} &= \begin{bmatrix} 0.8147 & 0.9134 & 0.2785 \\
0.9058 & 0.6324 & 0.5469 \\
0.1270 & 0.0975 & 0.9575 \end{bmatrix} \quad B_{\text{o}} = \begin{bmatrix} -0.6787 \\ -0.7577 \\ -0.7431 \end{bmatrix} \label{equ:44_1} \\[2mm]
C_1 &= \,\begin{bmatrix}\,1 & 0 & 0 \,\,\end{bmatrix} \quad\quad D_1 = \,\,0 \label{equ:44_2}  \\[2mm]
C_2 &= \begin{bmatrix} 0 & 1 & 0 \\
0 & 0 & 0 \end{bmatrix} \quad\quad D_2 = \begin{bmatrix} 0 \\
1 \end{bmatrix} \quad\quad y_{\max} = \begin{bmatrix} 1 \\ 0.5 \end{bmatrix} \label{equ:44_3} 
\end{align}
\end{subequations} \normalsize
The parameters in \eqref{equ:44_2} mean that we want to minimize the effect of disturbance input $w(t)$ on the first state $x_1(t)$, and the parameters in \eqref{equ:44_3} mean that the second state $x_2(t)$ and the control input $u(t)$ should satisfy the constraints $x_2(t) \le 1$ and $u(t) \le 0.5$. Assume $(A_{\text{o}},B_{\text{o}})$ is unknown. We simulate the system and collect $J = 100$ data of $(x^+_j,x_j,u_j)$ with random disturbance input $w_j$ that satisfies $\|w_j\| \le \varepsilon = 10^{-2}$. Then, we construct the matrices in Assumption~1 as $\Phi_{11} = (J \varepsilon^2) I$, $\Phi_{12} = 0$, and $\Phi_{22} = -I$.

We implement three data-driven approaches using the same set of collected data $\{(x^+_j,x_j,u_j)\}_{j = 1}^{J}$ to control the system: 1) the data-driven $\mathcal{H}_{\infty}$ control approach of \cite{van2020noisy} (which does not handle constraints), 2) the constrained $\mathcal{H}_{\infty}$ control based on \eqref{equ:20} (without moving-horizon implementation), and 3) the moving-horizon $\mathcal{H}_{\infty}$ control based on \eqref{equ:35}.

We consider the initial condition $x(0) = [0.95,0,0]^{\top}$. For~\eqref{equ:20} and \eqref{equ:35}, to illustrate the results of Theorems~1 and~2, we use the parameters $\sigma_0 = 10^{-2}$, $\sigma_{t+1} = \sigma_t - \|w(t)\|^2$, and $r_t = 10$ for all $t \in \mathbb{N}$. 

\begin{figure}[ht!]
\begin{center}
\begin{picture}(240.0, 148.0)
\put(  -2,  0){\epsfig{file=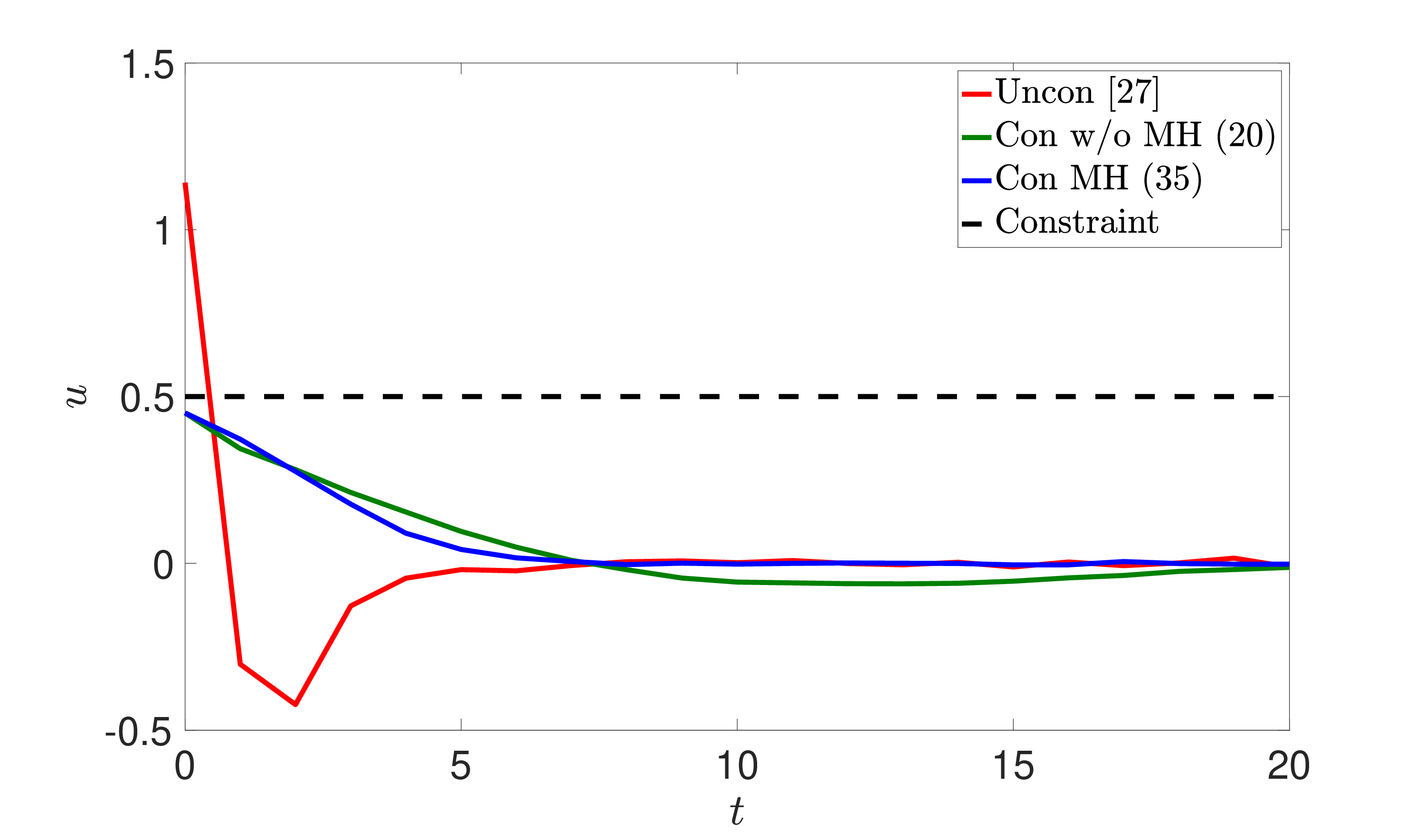,height=2in}}
\end{picture}
\end{center}
      \caption{Control input time history.}
      \label{fig:1}
\end{figure}

\begin{figure}[ht!]
\begin{center}
\begin{picture}(240.0, 148.0)
\put(  -2,  0){\epsfig{file=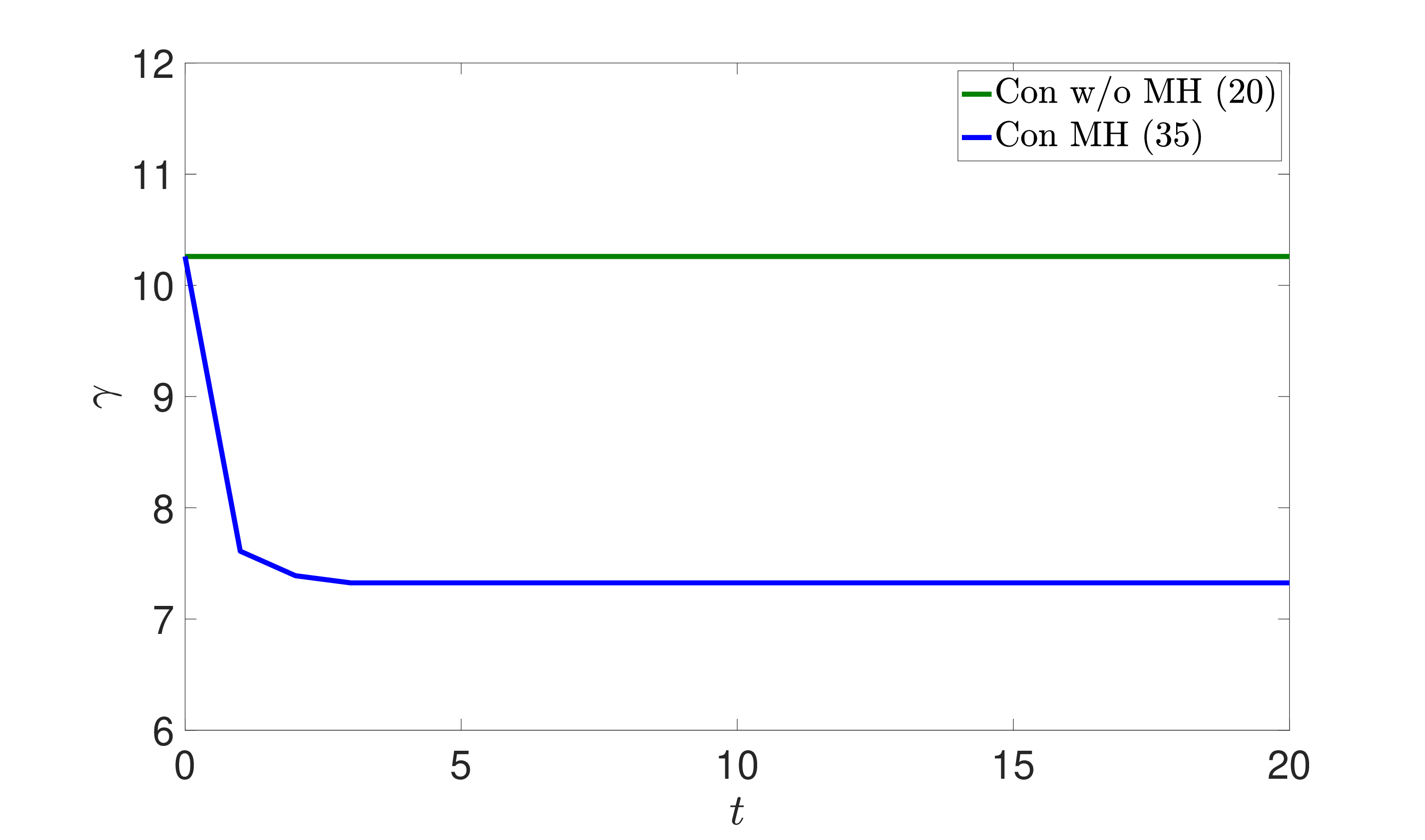,height=2in}}
\end{picture}
\end{center}
      \caption{Disturbance attenuation level time history.}
      \label{fig:2}
      \vspace{6mm}
\end{figure}

The system is stabilized and $x(t)$ converges to $0$ under each of the three approaches. However, as shown in Fig.~\ref{fig:1}, the control using the approach of \cite{van2020noisy} violates the constraint $u(t) \le 0.5$. This is expected because the approach of \cite{van2020noisy} does not consider any constraints. In contrast, \eqref{equ:20} and \eqref{equ:35} both satisfy the constraints, illustrating the effectiveness of our proposed approach for handling constraints. We then compare the time history of $\gamma$ using the constrained $\mathcal{H}_{\infty}$ control based on \eqref{equ:20} (without moving-horizon implementation) versus that using the moving-horizon $\mathcal{H}_{\infty}$ control based on \eqref{equ:35} in Fig.~\ref{fig:2}. Note that a smaller $\gamma$ indicates a stronger attenuation of the effect of disturbance input $w(t)$ on performance output $y_1(t)$. It can be seen that both approaches start with a large $\gamma$. This is because the control has to sacrifice some disturbance attenuation performance for satisfying the constraints on $x_2(t)$ and $u(t)$. Because the approach of \eqref{equ:20} uses a fixed gain, the disturbance attenuation performance level $\gamma$ remains constant over time. In contrast, as the state $x(t)$ and the control $u(t)$ are getting farther away from the constraint boundaries, the moving-horizon approach based on \eqref{equ:35} adjusts the gain and achieves a lower $\gamma$, illustrating the effectiveness of our proposed moving-horizon approach for improving performance while satisfying constraints.

%%%%%%%%%%%%%%%%%%%%%%%%%%%%%%%%%%%%%%%%%%%%%%%%%%%%%%%%%%%%%
%%%%%%%%%%%%%%%%%%%%%%%%%%%%%%%%%%%%%%%%%%%%%%%%%%%%%%%%%%%%%

%%%%%%%%%%%%%%%%%%%%%%%%%%%%%%%%%%%%%%%%%%%%%%%%%%%%%%%%%%%%%
%%%%%%%%%%%%%%%%%%%%%%%%%%%%%%%%%%%%%%%%%%%%%%%%%%%%%%%%%%%%%

In this example, the optimization problem \eqref{equ:35} is feasible at every time step, which is consistent with the recursive feasibility result of Lemma~4 and Theorem~2. The average computation time per step for solving \eqref{equ:35} using the MATLAB-based LMI solver {\sf mincx} on a MacBook~Air (M1~CPU, 8~GB~RAM) is {\sf 12.5~ms}, indicating the computational feasibility of the approach.

\section{Conclusions}\label{sec:6}

In this paper, we proposed a novel data-driven moving-horizon control approach for constrained systems. The~approach optimizes $\mathcal{H}_{\infty}$-type disturbance rejection while satisfying constraints. We established theoretical guarantees of the approach regarding closed-loop stability, disturbance attenuation, constraint satisfaction under noisy offline data, and online problem recursive feasibility. The effectiveness of the approach has been illustrated with a numerical example. Future work includes applying the proposed approach to practical control engineering problems.

%\section*{Acknowledgments}

\bibliographystyle{ieeetr}
\bibliography{ref}

\begin{thebibliography}{10}

\bibitem{zhou1996robust}
K.~Zhou, J.~C. Doyle, and K.~Glover, {\em Robust and optimal control}.
\newblock USA: Prentice Hall, 1996.

\bibitem{maciejowski2002predictive}
J.~M. Maciejowski, {\em Predictive control: with constraints}.
\newblock USA: Prentice Hall, 2002.

\bibitem{chen2004disturbance}
H.~Chen and C.~Scherer, ``Disturbance attenuation with actuator constraints by moving horizon {$H_{\infty}$} control,'' {\em IFAC Proceedings}, vol.~37, no.~1, pp.~415--420, 2004.

\bibitem{chen2006moving}
H.~Chen and C.~W. Scherer, ``Moving horizon {$H_{\infty}$} control with performance adaptation for constrained linear systems,'' {\em Automatica}, vol.~42, no.~6, pp.~1033--1040, 2006.

\bibitem{lee2008robust}
S.-M. Lee and J.~H. Park, ``Robust {$H_{\infty}$} model predictive control for uncertain systems using relaxation matrices,'' {\em International Journal of Control}, vol.~81, no.~4, pp.~641--650, 2008.

\bibitem{orukpe2011towards}
P.~E. Orukpe, ``Towards a less conservative model predictive control based on mixed {$H_2/H_{\infty}$} control approach,'' {\em International Journal of Control}, vol.~84, no.~5, pp.~998--1007, 2011.

\bibitem{huang2014mixed}
H.~Huang, D.~Li, and Y.~Xi, ``Mixed {$H_2/H_{\infty}$} robust model predictive control with saturated inputs,'' {\em International Journal of Systems Science}, vol.~45, no.~12, pp.~2565--2575, 2014.

\bibitem{benallouch2014h}
M.~Benallouch, G.~Schutz, D.~Fiorelli, and M.~Boutayeb, ``{$H_{\infty}$} model predictive control for discrete-time switched linear systems with application to drinking water supply network,'' {\em Journal of Process Control}, vol.~24, no.~6, pp.~924--938, 2014.

\bibitem{song2016mixed}
Y.~Song, X.~Fang, and Q.~Diao, ``Mixed {$H_2/H_{\infty}$} distributed robust model predictive control for polytopic uncertain systems subject to actuator saturation and missing measurements,'' {\em International Journal of Systems Science}, vol.~47, no.~4, pp.~777--790, 2016.

\bibitem{song2017robust}
Y.~Song, Z.~Wang, D.~Ding, and G.~Wei, ``Robust {$H_2/H_{\infty}$} model predictive control for linear systems with polytopic uncertainties under weighted {MEF-TOD} protocol,'' {\em IEEE Transactions on Systems, Man, and Cybernetics: Systems}, vol.~49, no.~7, pp.~1470--1481, 2017.

\bibitem{zhang2018robust}
Y.~Zhang, C.-C. Lim, and F.~Liu, ``Robust mixed {$H_2/H_{\infty}$} model predictive control for {Markov} jump systems with partially uncertain transition probabilities,'' {\em Journal of the Franklin Institute}, vol.~355, no.~8, pp.~3423--3437, 2018.

\bibitem{shokrollahi2021robust}
A.~Shokrollahi and S.~Shamaghdari, ``Robust {$H_{\infty}$} model predictive control for constrained {Lipschitz} non-linear systems,'' {\em Journal of Process Control}, vol.~104, pp.~101--111, 2021.

\bibitem{lewis2012reinforcement}
F.~L. Lewis, D.~Vrabie, and K.~G. Vamvoudakis, ``Reinforcement learning and feedback control: Using natural decision methods to design optimal adaptive controllers,'' {\em IEEE Control Systems Magazine}, vol.~32, no.~6, pp.~76--105, 2012.

\bibitem{recht2019tour}
B.~Recht, ``A tour of reinforcement learning: The view from continuous control,'' {\em Annual Review of Control, Robotics, and Autonomous Systems}, vol.~2, pp.~253--279, 2019.

\bibitem{rosolia2017learning}
U.~Rosolia and F.~Borrelli, ``Learning model predictive control for iterative tasks. {A} data-driven control framework,'' {\em IEEE Transactions on Automatic Control}, vol.~63, no.~7, pp.~1883--1896, 2017.

\bibitem{hewing2020learning}
L.~Hewing, K.~P. Wabersich, M.~Menner, and M.~N. Zeilinger, ``Learning-based model predictive control: Toward safe learning in control,'' {\em Annual Review of Control, Robotics, and Autonomous Systems}, vol.~3, pp.~269--296, 2020.

\bibitem{coulson2019data}
J.~Coulson, J.~Lygeros, and F.~D{\"o}rfler, ``Data-enabled predictive control: In the shallows of the {DeePC},'' in {\em 18th European Control Conference}, pp.~307--312, IEEE, 2019.

\bibitem{elokda2021data}
E.~Elokda, J.~Coulson, P.~N. Beuchat, J.~Lygeros, and F.~D{\"o}rfler, ``Data-enabled predictive control for quadcopters,'' {\em International Journal of Robust and Nonlinear Control}, vol.~31, no.~18, pp.~8916--8936, 2021.

\bibitem{huang2021decentralized}
L.~Huang, J.~Coulson, J.~Lygeros, and F.~D{\"o}rfler, ``Decentralized data-enabled predictive control for power system oscillation damping,'' {\em IEEE Transactions on Control Systems Technology}, vol.~30, no.~3, pp.~1065--1077, 2021.

\bibitem{chinde2022data}
V.~Chinde, Y.~Lin, and M.~J. Ellis, ``Data-enabled predictive control for building {HVAC} systems,'' {\em Journal of Dynamic Systems, Measurement, and Control}, vol.~144, no.~8, p.~081001, 2022.

\bibitem{li2023minimum}
N.~Li, E.~Taheri, I.~Kolmanovsky, and D.~Filev, ``Minimum-time trajectory optimization with data-based models: A linear programming approach,'' {\em arXiv preprint arXiv:2312.05724}, 2023.

\bibitem{baros2022online}
S.~Baros, C.-Y. Chang, G.~E. Colon-Reyes, and A.~Bernstein, ``Online data-enabled predictive control,'' {\em Automatica}, vol.~138, p.~109926, 2022.

\bibitem{dai2022cloud}
L.~Dai, T.~Huang, R.~Gao, Y.~Zhang, and Y.~Xia, ``Cloud-based computational data-enabled predictive control,'' {\em IEEE Internet of Things Journal}, vol.~9, no.~24, pp.~24949--24962, 2022.

\bibitem{berberich2020data}
J.~Berberich, J.~K{\"o}hler, M.~A. M{\"u}ller, and F.~Allg{\"o}wer, ``Data-driven model predictive control with stability and robustness guarantees,'' {\em IEEE Transactions on Automatic Control}, vol.~66, no.~4, pp.~1702--1717, 2020.

\bibitem{coulson2021distributionally}
J.~Coulson, J.~Lygeros, and F.~D{\"o}rfler, ``Distributionally robust chance constrained data-enabled predictive control,'' {\em IEEE Transactions on Automatic Control}, vol.~67, no.~7, pp.~3289--3304, 2021.

\bibitem{huang2023robust}
L.~Huang, J.~Zhen, J.~Lygeros, and F.~D{\"o}rfler, ``Robust data-enabled predictive control: Tractable formulations and performance guarantees,'' {\em IEEE Transactions on Automatic Control}, vol.~68, no.~5, pp.~3163--3170, 2023.

\bibitem{van2020noisy}
H.~J. van Waarde, M.~K. Camlibel, and M.~Mesbahi, ``From noisy data to feedback controllers: Nonconservative design via a matrix {S}-lemma,'' {\em IEEE Transactions on Automatic Control}, vol.~67, no.~1, pp.~162--175, 2020.

\bibitem{aastrom1977theory}
K.~J. {\AA}str{\"o}m, U.~Borisson, L.~Ljung, and B.~Wittenmark, ``Theory and applications of self-tuning regulators,'' {\em Automatica}, vol.~13, no.~5, pp.~457--476, 1977.

\bibitem{bristow2006survey}
D.~A. Bristow, M.~Tharayil, and A.~G. Alleyne, ``A survey of iterative learning control,'' {\em IEEE Control Systems Magazine}, vol.~26, no.~3, pp.~96--114, 2006.

\bibitem{hou2013model}
Z.-S. Hou and Z.~Wang, ``From model-based control to data-driven control: Survey, classification and perspective,'' {\em Information Sciences}, vol.~235, pp.~3--35, 2013.

\bibitem{de1992discrete}
C.~E. de~Souza and L.~Xie, ``On the discrete-time bounded real lemma with application in the characterization of static state feedback {$H_{\infty}$} controllers,'' {\em Systems \& Control Letters}, vol.~18, no.~1, pp.~61--71, 1992.

\bibitem{tuy1992nonconvex}
H.~Tuy, ``On nonconvex optimization problems with separated nonconvex variables,'' {\em Journal of Global Optimization}, vol.~2, pp.~133--144, 1992.

\end{thebibliography}

\end{document}